\documentclass[11pt,draft]{amsart}
\usepackage{amsfonts}
\usepackage{amsthm}
\usepackage{amsmath}
\usepackage{graphicx}
\usepackage{amscd,amssymb,amsthm}

\newcounter{minutes}
\divide\time by 60
\newcounter{hours}
\multiply\time by 60 \addtocounter{minutes}{-\time}

\newcommand{\dt}{{\rm d}t}
\newcommand{\dz}{{\rm d}z}

\newcommand{\doo}{{\rm d}}

\title{The generalized Marcum $Q-$function: an orthogonal
polynomial approach}

\author{Szil\'ard Andr\'as}
\address{Department of Applied Mathematics, Babe\c{s}-Bolyai
University, Cluj-Napoca 400084, Romania} \email{andraszk@yahoo.com}

\author[\'Arp\'ad Baricz]{\'Arp\'ad Baricz}
\address{Department of Economics, Babe\c{s}-Bolyai University,
Cluj-Napoca 400591, Romania} \email{bariczocsi@yahoo.com}

\author{Yin Sun}
\address{State Key Laboratory on Microwave and Digital Communications,
Tsinghua National Laboratory for Information Science and Technology
and Department of Electronic Engineering, Tsinghua University,
Beijing 100084, China} \email{sunyin02@mails.tsinghua.edu.cn}

\keywords{Marcum $Q-$function, Laguerre polynomials, modified Bessel
functions.}

\begin{document}

\begin{center}
\texttt{File:~\jobname .tex,
          printed: \number\year-0\number\month-0\number\day,
          \thehours.\ifnum\theminutes<10{0}\fi\theminutes}
\end{center}

\maketitle

\begin{abstract}
A novel power series representation of the generalized Marcum
$Q-$function of positive order involving generalized Laguerre
polynomials is presented. The absolute convergence of the proposed
power series expansion is showed, together with a convergence speed
analysis by means of truncation error. A brief review of related
studies and some numerical results are also provided.
\end{abstract}

\section{\bf Introduction}
\setcounter{equation}{0}

For $\nu$ real number let $I_{\nu}$ be denotes the modified Bessel
function \cite[p. 77]{watson} of the first kind of order $\nu,$
defined by
\begin{equation}\label{modbes}
I_{\nu}(t)=\sum_{n\geq0}\frac{(t/2)^{2n+\nu}}{n!\Gamma(\nu+n+1)},
\end{equation}
and let $b\mapsto Q_{\nu}(a,b)$ be the generalized Marcum
$Q-$function, defined by
\begin{equation}\label{maq}Q_{\nu}(a,b)=\frac{1}{a^{\nu-1}}\int_b^{\infty}t^{\nu}e^{-\frac{t^2+a^2}{2}}I_{\nu-1}(at)\dt,
\end{equation}
where $b\geq 0$ and $a,\nu>0.$ Here $\Gamma$ stands for the
well-known Euler gamma function. When $\nu=1,$ the function
$$b\mapsto Q_{1}(a,b)=\int_b^{\infty}te^{-\frac{t^2+a^2}{2}}I_{0}(at)\dt$$
is known in literature as the (first order) Marcum $Q-$function. The
Marcum $Q-$function and its generalization are frequently used in the
detection theories for radar systems \cite{marcum} and wireless communications \cite{Digham03,Digham07}, and have important applications in error
performance analysis of digital communication problems dealing with
partially coherent, differentially coherent, and non-coherent
detections \cite{sa,sun4}. Since, the precise
computations of the Marcum $Q-$function and generalized Marcum
$Q-$function are quite difficult, in the last few decades several
authors worked on precise and stable numerical calculation
algorithms for the functions. See the papers of Dillard \cite{dill},
Cantrell \cite{cantr}, Cantrell and Ojha \cite{Cantrell}, Shnidman
\cite{Shnidman}, Helstrom \cite{helstrom}, Temme \cite{temme1} and
the references therein. Moreover, many tight lower and upper bounds
for the Marcum $Q-$function and generalized Marcum $Q-$function were
proposed as simpler alternative evaluating methods or intermediate
results for further integrations. See, for example, the papers of
Simon \cite{simon}, Chiani \cite{chiani}, Simon and Alouini
\cite{sa1}, Annamalai and Tellambura \cite{ana}, Corazza and Ferrari
\cite{corazza}, Li and Kam \cite{lk}, Baricz \cite{baricz2}, Baricz
and Sun \cite{baricz1,baricz3}, Kapinas et al. \cite{mihos}, Sun et
al. \cite{sun2}, Li et al. \cite{li} and the references therein. In
this field, the order $\nu$ is usually the number of independent
samples of the output of a square-law detector, and hence in most of
the papers the authors deduce lower and upper bounds for the
generalized Marcum $Q-$function with order $\nu$ integer. On the other hand, based on the
papers \cite{Cantrell,marcum,Shnidman} there are introduced in the
Matlab 6.5 software the Marcum $Q-$function and positive integer order generalized
Marcum $Q-$function\footnote{See
\texttt{http://www.mathworks.com/access/helpdesk/help/toolbox/signal/marcumq.html}
for more details.}: \verb"marcumq(a,b)" computes the value of the
first order Marcum $Q-$function $Q_1(a,b)$ and \verb"marcumq(a,b,m)"
computes the value of the $m$th order generalized Marcum
$Q-$function $Q_{m}(a,b),$ defined by \eqref{maq}, where $m$ is a
positive integer. However, in some important applications, the order $\nu>0$ of the generalized Marcum $Q-$function is not necessarily an integer
number: The generalized Marcum $Q-$function is the
complementary cumulative distribution function or reliability
function of the non-central chi distribution with $2\nu$ degrees of
freedom \cite{cont,sun,sun2}. Moreover, real order generalized Marcum $Q-$function has been used to characterize small-scale channel fading distributions with line-of-sight channel components \cite{Loskot,Yacoub} or cross-channel correlations \cite{annama,anna,mihos,Khatalin,sa,Tellambura,Tellambura2}.

In this paper, we present a novel generalized Laguerre polynomial
series representation of the generalized Marcum $Q-$function, which
extends the result of the first order Marcum $Q-$function in Pent's
paper \cite{pent} to the case of the generalized Marcum $Q-$function
with real order $\nu>0$. We further show the absolute convergence of
the proposed power series expansion, together with a convergence
speed analysis by means of truncation error. A brief review of
related studies in the literature is provided, which may assist the
readers to get a more complete vision of this area. Finally, some
numerical results are provided as a complementary of these
theoretical analysis.
\section{\bf The generalized Marcum $Q-$function via Laguerre
polynomials} \setcounter{equation}{0}
\subsection{Novel series representation of the generalized Marcum
$Q-$function}
We start with the following well-known formula
\cite[p. 102]{szego}
\begin{equation}\label{sz1}
\sum_{n\geq0}\frac{L_n^{(\alpha)}(x)}{L_n^{(\alpha)}(0)}\frac{z^n}{n!}=\Gamma(\alpha+1)e^z(xz)^{-\frac{\alpha}{2}}J_{\alpha}(2\sqrt{xz}),
\end{equation}
where $x,z\in\mathbb{R}$ and $\alpha>-1.$ Here $J_{\alpha}$ stands
for the Bessel function of the first kind of order $\alpha,$
$L_n^{(\alpha)}$ is the generalized Laguerre polynomial of degree
$n$ and order $\alpha$, defined explicitly as
$$L_n^{(\alpha)}(x)=\frac{e^xx^{-\alpha}}{n!}\left(e^{-x}x^{n+\alpha}\right)^{(n)}=
\sum_{k=0}^n\frac{\Gamma(n+\alpha+1)}{\Gamma(k+\alpha+1)\Gamma(n-k+1)}\frac{(-x)^k}{k!}.$$
Changing in
\eqref{sz1} $z$ with $-z$ and taking into account the formula
$I_{\nu}(x)=\mathrm{i}^{-\nu}J_{\nu}(\mathrm{i}x)$ we obtain that
\cite{magnus}
\begin{equation}\label{t3}
\sum_{n\geq0}\frac{L_n^{(\alpha)}(x)}{L_n^{(\alpha)}(0)}\frac{(-1)^nz^n}{n!}=\Gamma(\alpha+1)e^{-z}(xz)^{-\frac{\alpha}{2}}I_{\alpha}(2\sqrt{xz}).
\end{equation}
Now, if we use
$$L_n^{(\alpha)}(0)=\frac{\Gamma(n+\alpha+1)}{\Gamma(\alpha+1)\Gamma(n+1)},$$
and replace $x$ with $a$ and $\alpha$ with $\nu-1,$ respectively, \eqref{t3} can be rewritten as
\begin{equation}\label{t1}
\left(\frac{z}{a}\right)^{\frac{\nu-1}{2}}e^{-z-a}I_{\nu-1}(2\sqrt{az})=
e^{-a}\sum_{n\geq0}(-1)^n\frac{L_n^{(\nu-1)}(a)}{\Gamma(\nu+n)}z^{n+\nu-1},
\end{equation}
which holds for all $a,\nu>0$ and $z\geq0.$

Now, consider the
following formula \cite{temme1,temme2}
\begin{equation}\label{t2}Q_{\nu}(\sqrt{2a},\sqrt{2b})=e^{-a}\sum_{n\geq0}\frac{\Gamma(\nu+n,b)}{\Gamma(\nu+n)}\frac{a^n}{n!}
=\int_b^{\infty}\left(\frac{z}{a}\right)^{\frac{\nu-1}{2}}e^{-z-a}I_{\nu-1}(2\sqrt{az})\dz,\end{equation}
where $a,\nu>0$ and $b\geq0.$ We note that the
function $b\mapsto Q_{\nu}(\sqrt{a},\sqrt{b}),$ defined by
$$Q_{\nu}(\sqrt{a},\sqrt{b})=\frac{1}{2}\int_b^{\infty}\left(\frac{z}{a}\right)^{\frac{\nu-1}{2}}e^{-\frac{z+a}{2}}I_{\nu-1}(\sqrt{az})\dz,$$
is in fact the survival function (or the complementary of the
cumulative distribution function with respect to unity) of the
non-central chi-square distribution with $2\nu$ degrees of freedom
and non-centrality parameter $a.$ With other words, for all
$a,\nu>0$ and $b\geq0$ we have
\begin{equation}\label{eq1}
Q_{\nu}(\sqrt{a},\sqrt{b})=1-\frac{1}{2}\int_0^b\left(
\frac{z}{a}\right)^{\frac{\nu-1}{2}}e^{-\frac{z+a}{2}}
I_{\nu-1}(\sqrt{az})\dz.
\end{equation}
See \cite{sun} for more details. Combining \eqref{t1} with \eqref{eq1}
we obtain
\begin{align*}
Q_{\nu}(\sqrt{2a},\sqrt{2b})&=
1-\int_0^b\left(\frac{z}{a}\right)^{\frac{\nu-1}{2}}e^{-z-a}I_{\nu-1}(2\sqrt{az})\dz\\&=
1-\int_0^be^{-a}\sum_{n\geq0}(-1)^n\frac{L_n^{(\nu-1)}(a)}{\Gamma(\nu+n)}z^{n+\nu-1}\dz\\&\overset{(a)}{=}
1-e^{-a}\sum_{n\geq0}(-1)^n\frac{L_n^{(\nu-1)}(a)}{\Gamma(\nu+n)}\int_0^bz^{n+\nu-1}\dz\\&=
1-\sum_{n\geq0}(-1)^ne^{-a}\frac{L_n^{(\nu-1)}(a)}{\Gamma(\nu+n+1)}b^{n+\nu},
\end{align*}
where in $(a)$ the integration and summation can be interchanged, because
the series on the right-hand side of \eqref{t1} is uniformly
convergent for $0\leq z\leq b$. For more details see the last paragraph of Section \ref{sec22}. After some simple manipulation, we obtain a new formula of the generalized Marcum $Q-$function, i.e.,
\begin{equation}\label{mn1}
Q_{\nu}(a,b)=1-\sum_{n\geq0}(-1)^ne^{-\frac{a^2}{2}}\frac{L_n^{(\nu-1)}\left(\frac{a^2}{2}\right)}{\Gamma(\nu+n+1)}\left(\frac{b^2}{2}\right)^{n+\nu},
\end{equation}
valid for all $a,\nu>0$ and $b\geq0.$

In order to simplify the numerical evaluation of the series \eqref{mn1}, we consider the expression
$$P_{\nu,n}(a,b)=\frac{b^nL_{n}^{(\nu-1)}(a)}{\Gamma(\nu+n+1)},$$
which satisfies the recurrence relation
\begin{align*}P_{\nu,n+1}(a,b)=&\frac{(2n+\nu-a)b}{(n+1)(\nu+n+1)}P_{\nu,n}(a,b)\\&-\frac{(n+\nu-1)b^2}{(n+1)(\nu+n)(\nu+n+1)}P_{\nu,n-1}(a,b)\end{align*}
for all $a,\nu>0,$ $b\geq 0$ and $n\in\{1,2,3,\dots\}$, with the
initial conditions
$$P_{\nu,0}(a,b)=\frac{1}{\Gamma(\nu+1)}\ \ \ \mbox{and}\ \ \ P_{\nu,1}(a,b)=\frac{(\nu-a)b}{\Gamma(\nu+2)}.$$
Here, the recurrence relation for $P_{\nu,n}(a,b)$ were obtained
from the recurrence relation \cite[p. 101]{szego}
$$(n+1)L_{n+1}^{(\alpha)}(x)=(2n+\alpha+1-x)L_n^{(\alpha)}(x)-(n+\alpha)L_{n-1}^{(\alpha)}(x)$$
and the initial conditions from
$$L_0^{(\alpha)}(x)=1\ \ \ \mbox{and}\ \ \ L_1^{(\alpha)}(x)=-x+\alpha+1.$$
With the
help of the expression $P_{\nu,n}(a,b)$, \eqref{mn1} can be easily rewritten as
\begin{equation}\label{mn2}
Q_{\nu}(a,b)=1-\sum_{n\geq0}e^{-\frac{a^2}{2}}\left(\frac{b^2}{2}\right)^{\nu}P_{\nu,n}\left(\frac{a^2}{2},-\frac{b^2}{2}\right).
\end{equation}
\subsection{Convergence analysis of the new series representation}\label{sec22} We note that for $a>0,$ $\nu\geq 1$ and $b\geq0$ the
absolute convergence of the series in \eqref{mn1} or \eqref{mn2} can
be shown easily by using the following inequalities
\begin{align*}
&\left|\sum_{n\geq0}(-1)^ne^{-\frac{a^2}{2}}\frac{L_n^{(\nu-1)}\left(\frac{a^2}{2}\right)}{\Gamma(\nu+n+1)}\left(\frac{b^2}{2}\right)^{n+\nu}\right|\\
&\leq
e^{-\frac{a^2}{2}}\sum_{n\geq0}\frac{1}{\Gamma(\nu+n+1)}\left(\frac{b^2}{2}\right)^{n+\nu}\left|L_n^{(\nu-1)}\left(\frac{a^2}{2}\right)\right|\\
& \leq
e^{-\frac{a^2}{2}}\sum_{n\geq0}\frac{1}{\Gamma(\nu+n+1)}\left(\frac{b^2}{2}\right)^{n+\nu}
\frac{\Gamma(\nu+n)}{n!\Gamma(\nu)}e^{\frac{a^2}{4}}\\
& \leq
e^{-\frac{a^2}{4}}\frac{1}{\Gamma(\nu)}\left(\frac{b^2}{2}\right)^{\nu-1}\sum_{n\geq0}\frac{1}{(n+1)!}\left(\frac{b^2}{2}\right)^{n+1}\\
&=e^{-\frac{a^2}{4}}\frac{1}{\Gamma(\nu)}\left(\frac{b^2}{2}\right)^{\nu-1}\left(e^{\frac{b^2}{2}}-1\right)
\end{align*}
or
\begin{align*}
&\left|\sum_{n\geq0}(-1)^ne^{-\frac{a^2}{2}}\frac{L_n^{(\nu-1)}\left(\frac{a^2}{2}\right)}{\Gamma(\nu+n+1)}\left(\frac{b^2}{2}\right)^{n+\nu}\right|\\
&\leq
e^{-\frac{a^2}{2}}\sum_{n\geq0}\frac{1}{\Gamma(\nu+n+1)}\left(\frac{b^2}{2}\right)^{n+\nu}\left|L_n^{(\nu-1)}\left(\frac{a^2}{2}\right)\right|\\
& \leq
e^{-\frac{a^2}{2}}\sum_{n\geq0}\frac{1}{\Gamma(\nu+n+1)}\left(\frac{b^2}{2}\right)^{n+\nu}
\frac{\Gamma(\nu+n)}{n!}e^{\frac{a^2}{4}}\left(\frac{a^2}{4}\right)^{1-\nu}\\
& \leq
e^{-\frac{a^2}{4}}\left(\frac{2b^2}{a^2}\right)^{\nu-1}\sum_{n\geq0}\frac{1}{(n+1)!}\left(\frac{b^2}{2}\right)^{n+1}\\
&=e^{-\frac{a^2}{4}}\left(\frac{2b^2}{a^2}\right)^{\nu-1}\left(e^{\frac{b^2}{2}}-1\right),
\end{align*}
which contain the known inequalities of Szeg\H o \cite{szego} for
generalized Laguerre polynomials
$$\left|L_n^{\alpha}(x)\right|\leq \frac{\Gamma(\alpha+n+1)}{n!\Gamma(\alpha+1)}e^{\frac{x}{2}}$$
and of Love \cite{love}
$$\left|L_n^{\alpha}(x)\right|\leq \frac{\Gamma(\alpha+n+1)}{n!}\left(\frac{x}{2}\right)^{-\alpha}e^{\frac{x}{2}},$$
where in both of the inequalities $\alpha\geq0,$ $x>0$ and
$n\in\{0,1,2,\dots\}.$

Moreover, for $a>0,$ $0<\nu\leq1$ and $b\geq0$ the absolute
convergence of the series in \eqref{mn1} or \eqref{mn2} can be shown
in a similar manner by using the following inequality
\begin{align*}
&\left|\sum_{n\geq0}(-1)^ne^{-\frac{a^2}{2}}\frac{L_n^{(\nu-1)}\left(\frac{a^2}{2}\right)}{\Gamma(\nu+n+1)}\left(\frac{b^2}{2}\right)^{n+\nu}\right|\\
&\leq
e^{-\frac{a^2}{2}}\sum_{n\geq0}\frac{1}{\Gamma(\nu+n+1)}\left(\frac{b^2}{2}\right)^{n+\nu}\left|L_n^{(\nu-1)}\left(\frac{a^2}{2}\right)\right|\\
& \leq
e^{-\frac{a^2}{2}}\sum_{n\geq0}\frac{1}{\Gamma(\nu+n+1)}\left(\frac{b^2}{2}\right)^{n+\nu}
\left(2-\frac{\Gamma(\nu+n)}{n!\Gamma(\nu)}\right)e^{\frac{a^2}{4}}\\
&=e^{-\frac{a^2}{4}}\sum_{n\geq0}\frac{1}{\nu+n}\left(\frac{2}{\Gamma(\nu+n)}-\frac{1}{n!\Gamma(\nu)}\right)\left(\frac{b^2}{2}\right)^{n+\nu}\\
& \leq e^{-\frac{a^2}{4}} \sum_{n\geq0}
\frac{2}{n!}\left(\frac{b^2}{2}\right)^{n+\nu}\\
&=
2e^{-\frac{a^2}{4}}\left(\frac{b^2}{2}\right)^{\nu}e^{\frac{b^2}{2}},
\end{align*}
which contains the classical inequality of Szeg\H o \cite{szego} for
generalized Laguerre polynomials
$$\left|L_n^{\alpha}(x)\right|\leq \left(2-\frac{\Gamma(\alpha+n+1)}{n!\Gamma(\alpha+1)}\right)e^{\frac{x}{2}}$$
where $-1<\alpha\leq0,$ $x>0$ and $n\in\{0,1,2,\dots\}.$ In addition
here we used the fact that for all fixed $n\in\{1,2,3,\dots\}$ the
function
$$\nu\mapsto\frac{1}{\nu+n}\left(\frac{2}{\Gamma(\nu+n)}-\frac{1}{n!\Gamma(\nu)}\right),$$
which maps $0$ into $2/n!,$ is decreasing on $(0,1]$ and
consequently for all $n\in\{0,1,2,\dots\}$ and $0<\nu\leq1$ we have
$$\frac{1}{\nu+n}\left(\frac{2}{\Gamma(\nu+n)}-\frac{1}{n!\Gamma(\nu)}\right)\leq \frac{2}{n!}.$$
We note that other sharper uniform bounds for generalized Laguerre
polynomials can be found in the papers of Love \cite{love},
Lewandowski and Szynal \cite{le}, Michalska and Szynal \cite{mi},
Pog\'any and Srivastava \cite{po}. See also the references therein.

Finally, note that by using the above uniform bounds for the
generalized Laguerre polynomials the uniform convergence of the
series on the right-hand side of \eqref{t1} can be shown easily for
$0\leq z\leq b.$ This is important because in order to obtain
\eqref{mn1} we have used tacitly that the series on the right-hand
side of \eqref{t1} is uniformly convergent and then we can
interchange the integration with summation. For example, if we use
the above Szeg\H o's uniform bound, then for all
$n\in\{0,1,2,\dots\},$ $a>0,$ $\nu\geq 1$ and $0\leq z\leq b$ we
have
$$\left|(-1)^n\frac{L_n^{(\nu-1)}(a)}{\Gamma(\nu+n)}z^{n}\right|\leq\frac{e^{\frac{a}{2}}}{\Gamma(\nu)}\frac{b^n}{n!}.$$
By the ratio test the series $e^b=\sum_{n\geq 0}b^n/n!$ is
convergent and thus in view of the Weierstrass M-test the original
series on the right-hand side of \eqref{t1} converges uniformly for
all $0\leq z\leq b.$

\subsection{Truncation error analysis}
For practical evaluations of our power series expansion, we need to approximate the generalized Marcum $Q-$function $Q_{\nu}(a,b)$ by the first
$n_0\in\{1,2,3,\dots\}$ terms of \eqref{mn1}, i.e.,
$$\hat{Q}_{\nu}(a,b)=1-\sum_{n=0}^{n_0}(-1)^ne^{-\frac{a^2}{2}}\frac{L_n^{(\nu-1)}\left(\frac{a^2}{2}\right)}
{\Gamma(\nu+n+1)}\left(\frac{b^2}{2}\right)^{n+\nu}.$$ We note that
the absolute value of the truncation error
$$\varepsilon_{\rm{t}}=Q_{\nu}(a,b)-\hat{Q}_{\nu}(a,b)=\sum_{n\geq n_0+1}(-1)^{n+1}e^{-\frac{a^2}{2}}\frac{L_n^{(\nu-1)}\left(\frac{a^2}{2}\right)}
{\Gamma(\nu+n+1)}\left(\frac{b^2}{2}\right)^{n+\nu}$$ can be upper
bounded by using the upper bounds for the generalized Laguerre
polynomials as in subsection 2.2. More precisely, by using the same
argument as in subsection 2.2 and Sewell's inequality \cite[p.
266]{mitri}
$$e^x-\sum_{k=0}^n\frac{x^k}{k!}\leq\frac{xe^x}{n}, \ \ \ n\in\{1,2,3,\dots\},\ x\geq0,$$ we can deduce
the following: if $a>0,$ $b\geq0$ and $\nu\geq 1,$ then
$$\left|\varepsilon_{\rm{t}}\right|\leq
e^{-\frac{a^2}{4}}\frac{1}{\Gamma(\nu)}\left(\frac{b^2}{2}\right)^{\nu-1}
\left[e^{\frac{b^2}{2}}-\sum_{n=0}^{n_0+1}\frac{1}{n!}\left(\frac{b^2}{2}\right)^{n}\right]
\leq\frac{e^{\frac{b^2}{2}-\frac{a^2}{4}}}{n_0+1}\frac{1}{\Gamma(\nu)}\left(\frac{b^2}{2}\right)^{\nu}$$
or
$$\left|\varepsilon_{\rm{t}}\right|\leq
e^{-\frac{a^2}{4}}\left(\frac{2b^2}{a^2}\right)^{\nu-1}
\left[e^{\frac{b^2}{2}}-\sum_{n=0}^{n_0+1}\frac{1}{n!}\left(\frac{b^2}{2}\right)^{n}\right]
\leq\frac{e^{\frac{b^2}{2}-\frac{a^2}{4}}}{n_0+1}\frac{b^2}{2}\left(\frac{2b^2}{a^2}\right)^{\nu-1}.$$
Similarly, it can be shown that if $a>0,$ $b\geq0$ and $0<\nu\leq1,$
then the absolute value of the truncation error is upper bounded as
follows
$$\left|\varepsilon_{\rm{t}}\right|\leq
2e^{-\frac{a^2}{4}}\left(\frac{b^2}{2}\right)^{\nu}\left[e^{\frac{b^2}{2}}-\sum_{n=0}^{n_0}
\frac{1}{n!}\left(\frac{b^2}{2}\right)^{n}\right]\leq
\frac{2e^{\frac{b^2}{2}-\frac{a^2}{4}}}{n_0}\left(\frac{b^2}{2}\right)^{\nu+1}.
$$
Observe that the above upper bounds of the absolute value of the
truncation error converge to zero at a speed of $1/n_0$. In
practice, we can use these upper bounds to decide the number of
terms, i.e. $n_0$, for achieving a pre-determined accuracy.
\subsection{A brief review of related studies}
As far as we know the formula \eqref{mn1}, or its equivalent form
\eqref{mn2}, is new. However, if we choose $\nu=1$ in \eqref{mn2},
then we reobtain the main result of Pent \cite{pent}
$$Q_{1}(a,b)=1-\frac{b^2}{2}\sum_{n\geq0}e^{-\frac{a^2}{2}}P_{n}\left(\frac{a^2}{2},-\frac{b^2}{2}\right),$$
where
$$P_{n}(a,b)=P_{1,n}(a,b)=\frac{b^nL_{n}(a)}{(n+1)!},$$
which for all $a>0,$ $b\geq 0$ and $n\in\{1,2,3,\dots\}$ satisfies
the recurrence relation
$$P_{n+1}(a,b)=\frac{(2n+1-a)b}{(n+1)(n+2)}P_{n}(a,b)-\frac{nb^2}{(n+1)^2(n+2)}P_{n-1}(a,b)$$
with the initial conditions
$$P_{0}(a,b)=1\ \ \ \mbox{and}\ \ \ P_{1}(a,b)=\frac{(1-a)b}{2}.$$ Here $L_n=L_n^{(0)}$ is the
classical Laguerre polynomial of degree $n.$

It should be mentioned here that another type of Laguerre expansions
for the Marcum $Q-$function was proposed in 1977 by Gideon and
Gurland \cite{Gideon}, which involves the lower incomplete gamma
function. This type of Laguerre expansions requires to use a
complementary result of \eqref{t3}, i.e.
\begin{equation}\label{t6}
\sum_{n\geq0}\frac{L_n^{(\alpha)}(z)}{L_n^{(\alpha)}(0)}\frac{(-1)^nx^n}{n!}
=\Gamma(\alpha+1)e^{-x}(xz)^{-\frac{\alpha}{2}}I_{\alpha}(2\sqrt{xz}).\end{equation}
Now by some simple manipulation we obtain
\begin{equation}\label{t4}
\left(\frac{z}{a}\right)^{\frac{\nu-1}{2}}e^{-z-a}I_{\nu-1}(2\sqrt{az})
=z^{\nu-1}e^{-z}\sum_{n\geq0}\frac{(-a)^n}{\Gamma(\nu+n)}L_n^{(\nu-1)}(z),
\end{equation}
which is equivalent to Tiku's result \cite{tiku}, available also as
equation (29.11) in the book \cite{cont}. By integrating \eqref{t4}
in $z$ and by using the differentiation formula \cite{magnus}
$$\frac{\doo}{\dz}\left[z^{\alpha+1}e^{-z}L_{n-1}^{(\alpha+1)}(z)\right]=nz^{\alpha}e^{-z}L_{n}^{(\alpha)}(z),$$
where $n\in\{1,2,3,\dots\},$ $\alpha>-1$ and $z\in\mathbb{R},$ we
can obtain another generalized Laguerre polynomial series expansion
of the generalized Marcum $Q-$function
$$Q_{\nu}(\sqrt{2a},\sqrt{2b})=1-\frac{1}{\Gamma(\nu)}\gamma(\nu,b)-\sum_{n\geq1}(-1)^ne^{-b}\frac{b^{\nu}L_{n-1}^{(\nu)}(b)}{n\Gamma(\nu+n)}a^n,$$
which in turn implies that
\begin{eqnarray}\label{t5}
Q_\nu(a,b)
\!\!\!\!\!\!\!\!\!\!\!&&=1-\frac{1}{\Gamma(\nu)}\gamma\left(\nu,\frac{b^2}{2}\right)-\sum_{n\geq
1}(-1)^ne^{-\frac{b^2}{2}}\left(\frac{b^2}{2}\right)^{\nu}\frac{L_{n-1}^{(\nu)}\left(\frac{b^2}{2}\right)}{n\Gamma(\nu+n)}\left(\frac{a^2}{2}\right)^n\nonumber\\
&&=\frac{1}{\Gamma(\nu)}\Gamma\left(\nu,\frac{b^2}{2}\right)-\sum_{n\geq
1}(-1)^ne^{-\frac{b^2}{2}}\left(\frac{b^2}{2}\right)^{\nu}\frac{L_{n-1}^{(\nu)}\left(\frac{b^2}{2}\right)}{n\Gamma(\nu+n)}\left(\frac{a^2}{2}\right)^n\nonumber\\
&&=\lim_{a\to0}Q_\nu(a,b)-\sum_{n\geq
1}(-1)^ne^{-\frac{b^2}{2}}\left(\frac{b^2}{2}\right)^{\nu}\frac{L_{n-1}^{(\nu)}\left(\frac{b^2}{2}\right)}{n\Gamma(\nu+n)}\left(\frac{a^2}{2}\right)^n,
\end{eqnarray}
where $\gamma(\cdot,\cdot)$ is the lower incomplete gamma function,
defined by
$$\gamma(a,x)=\int_0^xt^{a-1}e^{-t}\dt.$$
Here we used that
\begin{equation}\label{t7}
\Gamma(a,x)=\Gamma(a)-\gamma(a,x),
\end{equation}
and
$$\lim_{a\to0}Q_\nu(a,b)=\frac{1}{\Gamma(\nu)}\Gamma\left(\nu,\frac{b^2}{2}\right).$$
Some other Laguerre expansions for the Marcum $Q-$function are
provided in Gideon and Gurland's paper \cite{Gideon}, available also
as equation (29.13) of \cite{cont}. Moreover, a new unified Laguerre
polynomial-series-based distribution of small-scale fading envelope
and power was proposed recently by Chai and Tjhung \cite{chai},
which covers a wide range of small-scale fading distributions in
wireless communications. Many known Laguerre polynomial-series-based
probability density functions and cumulative distribution functions
of small-scale fading distributions are provided, which include the
multiple-waves-plus-diffuse-power fading, non-central chi and
chi-square, Nakagami-$m$, Rician (Nakagami-$n$), Nakagami-$q$
(Hoyt), Rayleigh, Weibull, Stacy, gamma, Erlang and exponential
distributions as special cases. See also \cite{sun3}, which contains
some corrections of formulas deduced in \cite{chai}. In particular,
\eqref{t5} is a special case of the unified cumulative distribution
function given in corrected form in \cite{sun3}. We note that the
expression of \eqref{t5} and the unified cumulative distribution
function in \cite{sun3} are quite different from our main result
\eqref{mn1} or \eqref{mn2}. This is because they are based on two
different Laguerre polynomial expansions of the modified Bessel
function of the first kind $I_\nu$ given in \eqref{t3} and
\eqref{t6}. Therefore, these Laguerre polynomial expansions are
expanded over different variables of the generalized Marcum
$Q-$function. Finally, we note that since Nakagami's work
\cite{naka} the Laguerre polynomial series expansions of various
probability density functions have been derived. We refer to the
papers of Esposito and Wilson \cite{esp}, Yu et al. \cite{yu}, Chai
and Tjhung \cite{chai} and to the references therein.

Finally, by using the infinite series representation of the modified Bessel
function of the first kind \eqref{modbes} and the formula
$$\int_{\alpha}^{\infty}t^me^{-\frac{t^2}{2}}\dt=2^{\frac{m-1}{2}}\Gamma\left(\frac{m+1}{2},\frac{\alpha^2}{2}\right),$$
where $\Gamma(\cdot,\cdot)$ is the upper incomplete gamma function,
defined by
$$\Gamma(a,x)=\int_x^{\infty}t^{a-1}e^{-t}\dt,$$
we easily obtain that
\begin{align}\label{maq2}
\nonumber
Q_{\nu}(a,b)&=\frac{1}{a^{\nu-1}}\int_b^{\infty}t^{\nu}e^{-\frac{t^2+a^2}{2}}\sum_{n\geq0}\frac{(at)^{2n+\nu-1}}{2^{2n+\nu-1}n!\Gamma(\nu+n)}\dt\\\nonumber
&=e^{-\frac{a^2}{2}}\sum_{n\geq0}\frac{a^{2n}}{2^{2n+\nu-1}n!\Gamma(\nu+n)}\int_b^{\infty}e^{-\frac{t^2}{2}}t^{2n+\nu-1}\dt\\\nonumber
&=e^{-\frac{a^2}{2}}\sum_{n\geq0}\frac{1}{n!}\left(\frac{a^2}{2}\right)^n\frac{\Gamma\left(\nu+n,\frac{b^2}{2}\right)}{\Gamma(\nu+n)}\\
&=1-\sum_{n\geq0}e^{-\frac{a^2}{2}}\left(\frac{a^2}{2}\right)^{n}\frac{\gamma\left(\nu+n,\frac{b^2}{2}\right)}{\Gamma(\nu+n)}.
\end{align}
We note that \eqref{maq2} is usually called the canonical
representation of the $\nu$th order generalized Marcum $Q-$function.
Recently, Annamalai and Tellambura \cite{annama} (see also
\cite{anna}) claimed that the series representation \eqref{maq2} is
new, however it appears already in 1993 in the paper of Temme
\cite{temme1}. See also Temme's book \cite{temme2} and Patnaik's
\cite{patnaik} result from 1949, which can be found also as equation
(29.2) in the book \cite{cont}. Interestingly, our novel series
representation \eqref{mn2} for the generalized Marcum $Q-$function
resembles to the series representation \eqref{maq2}.
\subsection{Numerical results} We now consider some numerical aspects of our
generalized Laguerre polynomial expansions \eqref{mn1} or
\eqref{mn2}. In practice, we usually need
to compute the detection probability for different values of $b$
with fixed $a$ to decide a proper detection threshold.
Since the generalized Laguerre polynomial in \eqref{mn1} is determined by only $a$, we can save computation time by storing the
values of the generalized Laguerre polynomials for computing the
generalized Marcum $Q-$function with different values of $b$.

The following tables contain some values of the generalized Marcum
$Q-$function calculated using \eqref{mn2} and using the Matlab
\verb"marcumq" function. For the considered choices of $a$ and $b$,
the numerical value of \eqref{mn2} is exactly the same with that of
the Matlab \verb"marcumq" function, if $\nu=1,3,5$ is integer. When
$\nu = 7.7$, the Matlab \verb"marcumq" function does not work, and
the numerical value of \eqref{mn2} is provided in the tables.
Finally, we note that more accurate intermediate terms are required
for larger $a$ and $b$.

{\footnotesize
\begin{center}
\begin{tabular}{|c|cccc|}
\hline
 $a=0.2, b=0.6$& $\nu=1$ & $\nu=3$ & $\nu=5$ & $\nu=7.7$\\
 \hline (\ref{mn2})&0.838249985438908 &0.999166310455636 &0.999998670306184 &  0.999999999927717\\
\hline marcumq &0.838249985438908 & 0.999166310455636 &0.999998670306184 &--- \\
\hline
\end{tabular}
\end{center}

\begin{center}
\begin{tabular}{|c|cccc|}
\hline
 $a=1.2, b=1.6$& $\nu=1$ & $\nu=3$ & $\nu=5$ & $\nu=7.7$\\
 \hline (\ref{mn2})&0.501536568390858 &0.916936068900377 & 0.994346394491553 &  0.999944937223540\\
\hline marcumq &0.501536568390858 & 0.916936068900377 & 0.994346394491553 &  --- \\
\hline
\end{tabular}
\end{center}

\begin{center}
\begin{tabular}{|c|cccc|}
\hline
 $a=2.2, b=2.6$& $\nu=1$ & $\nu=3$ & $\nu=5$ & $\nu=7.7$\\
 \hline (\ref{mn2})&0.426794627821735 &0.746459898209090 & 0.929671935077756 &  0.993735633182201\\
\hline marcumq &0.426794627821735 & 0.746459898209090 & 0.929671935077756 &  --- \\
\hline
\end{tabular}
\end{center}}

\subsection*{Acknowledgments} The work of S. Andr\'as was partially
supported by the Hungarian University Federation of Cluj. The
research of \'A. Baricz was supported by the J\'anos Bolyai Research
Scholarship of the Hungarian Academy of Sciences and by the Romanian
National Authority for Scientific Research CNCSIS-UEFISCSU, project
number PN-II-RU-PD\underline{ }388/2011. The work of Y. Sun was
supported by National Basic Research Program of China
(2007CB310608), National Natural Science Foundation of China
(60832008) and Lab project from Tsinghua National Lab on Information
Science and Technology (sub-project): Key technique for new
distributed wireless communications system.

\end{document}